\documentclass{amsart}
\usepackage{amsmath,amssymb,amsxtra,amsthm,amscd}
\usepackage[mathscr]{eucal}
\usepackage{bbm}
\usepackage{graphicx}
\usepackage[all,cmtip]{xy}
\usepackage{lmodern}
\usepackage{microtype}
\newif\ifShowLabels
\ShowLabelstrue
\newcommand{\TeXref}[1]{
\marginpar{\scriptsize \texttt{#1}}}
\ShowLabelsfalse

\DeclareMathOperator{\B}{\mathbf{B}}

\DeclareMathOperator{\im}{im}

\DeclareMathOperator{\LI}{\mathbf{LI}}

\DeclareMathOperator{\Mod}{\mathbf{Mod}}

\DeclareMathOperator*{\one}{1}
\newcommand{\onehatplace}[1]
{ \one^{\substack{#1 \\ \frown}} }

\DeclareMathOperator*{\bones}{\times}
\newcommand{\undertimes}[1]
{ \bones_{#1} }

\DeclareMathOperator*{\bowl}{\cup}
\newcommand{\undercup}[1]
{ \bowl_{#1} }

\DeclareMathOperator*{\arch}{\cap}
\newcommand{\undercap}[1]
{ \arch_{#1} }

\newcommand{\pull}
{\!\!\! -\!\!\! -\!\!\! -\!\!\!}

\DeclareMathOperator*{\holimprep}{holim}                       
\newcommand{\holim}[1]%
{\displaystyle\holimprep_{\substack{\leftarrow \pull - \\ #1}} \, }

\DeclareMathOperator*{\hocolimprep}{hocolim}                   
\newcommand{\hocolim}[1]%
{\displaystyle\hocolimprep_{\substack{- \pull \rightarrow \\ #1}} \, }

\DeclareMathOperator*{\plainlim}{lim}                           
\newcommand{\contralim}[1]%
{\displaystyle\plainlim_{\substack{\leftarrow \pull - \\ #1}} \, }

\DeclareMathOperator*{\plaincolim}{colim}                       
\newcommand{\colim}[1]%
{\displaystyle\plaincolim_{\substack{- \pull \rightarrow \\ #1}} \, }

\DeclareMathOperator*{\laxlimplain}{laxlim}                     
\newcommand{\laxlim}[1]%
{\displaystyle\laxlimplain_{\substack{\leftarrow \pull - \\ #1}} \, }

\providecommand{\bysame}{\makebox[3em]{\hrulefill}\thinspace}

\theoremstyle{plain}
\newtheorem{Thm}{Theorem}[section]

\newtheorem{Cor}[Thm]{Corollary}

\newtheorem{Lem}[Thm]{Lemma}
\newtheorem{Prop}[Thm]{Proposition}

\theoremstyle{definition}
\newtheorem{Def}[Thm]{Definition}

\newtheorem{Ex}[Thm]{Example}

\newtheorem{Rem}[Thm]{Remark}

\theoremstyle{remark}
\newtheorem{Not}[Thm]{Notation}

\newtheoremstyle{freestylethm}{6pt}{6pt}{\itshape}{}%
                {\bfseries}{}{.5em}{\thmnote{#3}}
\theoremstyle{freestylethm}

\newcommand{\refT}[1]{\textup{\ref{T:#1}}}
\newcommand{\refL}[1]{\textup{\ref{L:#1}}}
\newcommand{\refD}[1]{\textup{\ref{D:#1}}}
\newcommand{\refC}[1]{\textup{\ref{C:#1}}}

\newcommand{\refR}[1]{\textup{\ref{R:#1}}}

\newenvironment{ThmRef}[1]%
{ \begin{Thm} \label{T:#1}
\ifShowLabels \TeXref{T:#1} \fi }%
{ \end{Thm} }
\newenvironment{DefRef}[1]%
{ \begin{Def} \label{D:#1}
\ifShowLabels \TeXref{D:#1} \fi }%
{ \end{Def} }
\newenvironment{LemRef}[1]%
{ \begin{Lem} \label{L:#1}
\ifShowLabels \TeXref{L:#1} \fi }%
{ \end{Lem} }
\newenvironment{CorRef}[1]%
{ \begin{Cor} \label{C:#1}
\ifShowLabels \TeXref{C:#1} \fi }%
{ \end{Cor} }
\newenvironment{RemRef}[1]%
{ \begin{Rem} \label{R:#1}
\ifShowLabels \TeXref{R:#1} \fi }%
{ \end{Rem} }
{ \begin{Prop} \label{P:#1}
\ifShowLabels \TeXref{P:#1} \fi }%
{ \end{Prop} }
{ \begin{Ex} \label{E:#1}
\ifShowLabels \TeXref{E:#1} \fi  }%
{ \end{Ex} }
\newenvironment{NotRef}[1]%
{ \begin{Not} \label{N:#1}
\ifShowLabels \TeXref{N:#1} \fi }%
{ \end{Not} }

\newenvironment{ThmRefName}[2]%
{ \begin{Thm} [#2]\label{T:#1}
\ifShowLabels \TeXref{T:#1} \fi }%
{ \end{Thm} }
\newenvironment{DefRefName}[2]%
{ \begin{Def} [#2]\label{D:#1}
\ifShowLabels \TeXref{D:#1} \fi }%
{ \end{Def} }
{ \begin{Lem} [#2]\label{L:#1}
\ifShowLabels \TeXref{L:#1} \fi }%
{ \end{Lem} }
{ \begin{Cor} [#2]\label{C:#1}
\ifShowLabels \TeXref{C:#1} \fi }
{ \end{Cor} }
{ \begin{Rem} [#2]\label{R:#1}
\ifShowLabels \TeXref{R:#1} \fi }%
{ \end{Rem} }
{ \begin{Prop} [#2]\label{P:#1}
\ifShowLabels \TeXref{P:#1} \fi }%
{ \end{Prop} }
{ \begin{Ex} [#2]\label{E:#1}
\ifShowLabels \TeXref{E:#1} \fi }%
{ \end{Ex} }

\begin{document}
\title[Coarse coherence of metric spaces]{Coarse coherence of metric spaces and groups and its permanence properties}
\author[Boris Goldfarb]{Boris Goldfarb}
\address{Department of Mathematics and Statistics\\ SUNY\\ Albany\\ NY 12222}
\email{bgoldfarb@albany.edu}
\author[Jonathan L. Grossman]{Jonathan L. Grossman}
\address{Department of Mathematics and Statistics\\ SUNY\\ Albany\\ NY 12222}
\email{algrossman@albany.edu}
\date{\today}

\begin{abstract}
We introduce properties of metric spaces and, specifically, finitely generated groups with word metrics which we call \textit{coarse coherence} and \textit{coarse regular coherence}.
They are geometric counterparts of the classical algebraic notion of coherence and the regular coherence property of groups defined and studied by F. Waldhausen. 
The new properties can be defined in the general context of coarse metric geometry and are coarse invariants. In particular, they are quasi-isometry invariants of spaces and groups.  

We show that coarse regular coherence implies weak regular coherence, a weakening of regular coherence by G.~Carlsson and the first author.  The latter was introduced with the same goal as Waldhausen's, in order to perform computations of algebraic $K$-theory of group rings.
However, all groups known to be weakly regular coherent are also coarsely regular coherent. The new framework allows us to prove structural results by developing permanence properties, including the particularly important fibering permanence property, for coarse regular coherence.
\end{abstract}

\maketitle


\section{Precursors of coarse coherence}

Let $A$ be an associative ring.  All modules over $A$ we consider are left $A$-modules.  

\begin{DefRef}{HUYIV}
A \textit{presentation} of an $A$-module $E$ is an exact sequence  $F_2 \to F_1 \to E \to 0$ with both $F_1$ and $F_2$ free $A$-modules. It is a \textit{finite presentation} if the free modules are finitely generated.  More generally, one has the notion of a \textit{projective resolution} of $E$ which is an exact sequence 
\[
\ldots \longrightarrow P_n \longrightarrow \ldots \longrightarrow P_2 \longrightarrow P_1 \longrightarrow E \longrightarrow 0
\]
where all $P_i$ are projective $A$-modules.  The projective resolution is of \textit{finite type} if the projective modules are finitely generated.  It is called \textit{finite} if there is a number $n$ such  that the modules $P_i =0$ for $i > n$.

A module is said to be of type $\mathrm{FP}_{\infty}$ or have \textit{finite projective dimension} if
it has a projective resolution of finite type.
The ring $A$ has \textit{finite global dimension}
if there is a number $n$ such that every finitely generated $A$-module has a finite projective resolution of length $n$. 

The ring $A$ is called \textit{coherent} if every finitely presented $A$-module is of type $\mathrm{FP}_{\infty}$.
A coherent ring $A$ is called \textit{regular coherent} if each projective resolution of finite type over $A$ is chain homotopy equivalent
to a finite projective resolution.
Restricting further, a \textit{regular Noetherian} ring is a regular coherent ring which is Noetherian in the usual sense---a submodule of any finitely generated module over $A$ is finitely generated.

The most important regular Noetherian ring for applications in geometric topology is the ring of integers $\mathbb{Z}$.
\end{DefRef}

\begin{DefRefName}{WaldRC}{Waldhausen \cite{fW:78}}
	A group $\Gamma$ is \textit{regular coherent} if the group algebra $R[\Gamma]$ is
regular coherent for any choice of a regular Noetherian ring $R$.
\end{DefRefName}

The collection $\frak{X}$ of regular coherent groups includes free groups, free abelian groups,
torsion-free one relator groups, fundamental groups of
submanifolds of the three-dimensional sphere, and their various
amalgamated products and HNN extensions and so, in particular, the
fundamental groups of submanifolds of the three-dimensional sphere.  Waldhausen used this property to compute the algebraic $K$-theory of
regular coherent groups. 

Two remarks regarding Waldhausen's regular coherence are in order.

(1) The regular coherence property seems to be very special: simply constructing
individual non-projective finite dimensional modules over group rings is hard.  

(2) The collection $\frak{X}$ is not well-understood structurally beyond the portion identified by Waldhausen.  For example, it is unknown whether $\frak{X}$ is closed under products.  While all groups in $\frak{X}$ are necessarily torsion-free, it is unknown if there is a torsion-free group outside of it.

Waldhausen asked if a weaker property of the group
ring would suffice in his argument,
see for example the paragraph after the proof of Theorem 11.2 in~\cite{fW:78}.
This paper is a response to that question.

\medskip

A weakening of regular coherence (\textit{weak regular coherence}) was introduced in \cite{gCbG:04a,gCbG:04,gCbG:16,bG:13} with essentially the same goal as Waldhausen's which was a computation of the $K$-theory of weakly regular coherent groups.  However, in contrast with the situation for $\frak{X}$, the class of weakly regular coherent groups is known to be very large.  It includes groups that admit a finite classifying space and have straight finite decomposition complexity and so, in particular, groups that have finite asymptotic dimension.
The notion of weak coherence is more technical to define; we do that as part of the narrative in section~4.
Unfortunately one realizes quickly that even though the notion of weak regular coherence can be verified for a large family of groups, it is not amenable to proving structural permanence results.  

\medskip

\textbf{The goal and structure of the paper.}
Our goal is to define a genuinely coarse geometric property of metric spaces that ensures that finitely generated groups with word metrics that have this property also have the weak regular coherence property.  Then we show that the class of coarsely regular coherent groups $\frak{Y}$ is closed under many natural geometric operations. This collection of results has became known as \textit{permanence properties} \cite{eG:14} in the literature.  For example, invariance of $\frak{Y}$ under finite products is a simple consequence of this paper.  

We start by defining the general non-equivariant property of metric spaces called \textit{coarse coherence} in section 2 and stating main technical results.
Section 3 contains proofs of permanence results for coarse coherence.  
In the last section 4 of the paper, we define \textit{coarse regular coherence} which is a property of groups.   We relate this property to weak regular coherence and finally state and prove theorems about coarse regular coherence together with some immediate applications. 
 
\medskip

\textbf{Acknowledgements.}
We would like to thank Daniel Kasprowski, Marco Varisco, and the referee for comments that improved the strength of the results and the precision of the narrative.

\section{Definition of coarse coherence}

Let $X$ be a metric space, $R$ be a ring.    

\begin{DefRef}{Basic}
An $X$-filtered $R$-module is a covariant functor $F \colon P(X) \to \Mod_R$  from the power set of $X$ ordered by inclusion to the category of $R$-modules and injective homomorphisms. It will be convenient to view $F$ as the value $F(X)$ filtered by submodules associated to subsets $S \subset X$. We will always assume that the value of $F$ on the empty subset is $0$. 

The notation $S[r]$ stands for the
metric $r$-enlargement of $S$ in $X$.
So, in particular, $x[r]$ is the closed metric ball of radius $r$ centered at $x$.  

\begin{enumerate}
\item $F$ is called \textit{lean} or $D$-\textit{lean} if there is a number $D \ge 0$ such that
\[
F(S) \subset \sum_{x \in S} F(x[D])
\]
for every subset $S$ of $X$.
\item $F$ is called \textit{scattered} or $\delta$-\textit{scattered} if there is a number $\delta \ge 0$ such that
\[
F(X) \subset \sum_{x \in X} F(x[\delta]).
\]
\item $F$ is called \textit{insular} or $d$-\textit{insular} if there is a
number $d \ge 0$ such that
\[
F(S) \cap F(U) \subset F(S[d] \cap U[d])
\]
for every pair of subsets $S$, $U$ of $X$.
\item $F$ is \textit{locally finitely generated} if $F (S)$ is a finitely generated $R$-module for
every bounded subset $S \subset X$.
\end{enumerate}
\end{DefRef}   

\begin{RemRef}{LvS}
We note that being scattered is a consequence of being lean but not the other way around.
\end{RemRef}

\begin{DefRef}{UIOL}
An $R$-homomorphism $f \colon F \to F'$ of $X$-filtered modules is \textit{controlled} if there is a fixed number $b \ge 0$ such that the
image $f (F (S))$ is a submodule of $F' (S [b])$
for all subsets $S$ of $X$.
\end{DefRef}

The filtered modules that are lean and insular form a category $\LI (X,R)$, where the morphisms are the controlled $R$-homomorphisms.
The category $\B (X,R)$ is the full
subcategory of $\LI (X,R)$ on the locally finitely generated objects.
Basic properties of this category can be found in section 3.1 of \cite{gCbG:15}.  

Let us point out that the geometric features we define are independent of the assumption that $R$ is Noetherian which is used throughout \cite{gCbG:15}.

\begin{Ex}
	Given any ring $R$ and any metric space $X$, one can consider a special kind of filtrations which go back to the original geometric modules of Pedersen-Weibel.
	The \textit{geometric modules} are a collection of choices $F_x$ which are free finitely generated $R$-modules associated to each point $x$ in $X$, with the requirement that only finitely many $F_x$ are non-zero for $x$ from a bounded subset $S$.	Now the module $F (X) = \bigoplus F_x$ is assigned the filtration given by $F(S)=\bigoplus_{x \in S} F_x$.
	Let us denote the full subcategory of geometric modules in $\B (X,R)$ by $\mathcal{B} (X,R)$.
	\end{Ex}
	
	Within $\LI (X,R)$ there is a special kind of morphisms which were used in the description of the exact structure in $\LI (X,R)$. 

\begin{DefRef}{bic}
	A homomorphism $f$ is \textit{bicontrolled} if, for some fixed $b \ge 0$, in addition to inclusions of submodules
$f (F (S)) \subset F' (S [b])$,
there are inclusions
$f (F)
\cap F' (S) \subset  f F (S[b])$
for all subsets $S \subset X$.
\end{DefRef}

Here is a list of facts about lean and insular modules.

\begin{ThmRef}{lninpres}
Let
\[
0 \to E' \xrightarrow{\ f \ } E \xrightarrow{\ g \ } E'' \to 0
\]
be an exact sequence of $X$-filtered $R$-modules where $f$ and $g$ are bicontrolled.

\begin{enumerate}
\item If the object $E$ is
lean then $E''$ is lean.

\item If $E$ is
insular then $E'$ is insular.

\item If $E$ is
insular and $E'$ is lean then $E''$ is insular.

\item If both $E'$ and $E''$ are lean then $E$ is lean.
 
\item If both $E'$ and $E''$ are insular then $E$ is insular.
\end{enumerate}
\end{ThmRef}

\begin{proof}
	This is a summary of results from section 3.1 in \cite{gCbG:15}.
\end{proof}

There are several viable conditions on $X$ that enforce a version of the ``missing'' item in this theorem.  All of them can be viewed as relaxations of the algebraic coherence property when $X$ is a group $\Gamma$ with a word metric.
In this case the filtered $R$-modules are $\Gamma$-equivariant, thus becoming $R[\Gamma]$-modules, with the maps $f$ and $g$ being $R[\Gamma]$-homomorphisms.

The first condition is the ``missing'' item itself.

\begin{DefRefName}{Coherence}{Coherence of Metric Spaces}
A metric space $X$ is \textit{coherent} if in any exact sequence 
\[
0 \to E' \xrightarrow{\ f \ } E \xrightarrow{\ g \ } E'' \to 0
\]
of $X$-filtered $R$-modules where $f$ and $g$ are both bicontrolled maps, the combination of $E$ being lean and $E''$ being insular implies that $E'$ is necessarily lean.
\end{DefRefName}

For example, it is shown in Grossman \cite{jG:18} that the real line with the standard metric is a coherent metric space.
It is likely that all groups from Waldhausen's class $\mathfrak{X}$ are coherent.
Loc. cit. also contains basic properties of coherence such as coarse, and therefore quasi-isometry, invariance.  
It turns out however that for the most desired permanence properties this notion is too restrictive.

The following definition isolates the most important property in terms of its algebraic impact which also happens to be amenable to fibering permanence results.  

\begin{DefRefName}{WeakCoherence}{Coarse Coherence}
A metric space $X$ is \textit{coarsely coherent} if in any exact sequence 
\[
0 \to E' \xrightarrow{\ f \ } E \xrightarrow{\ g \ } E'' \to 0
\]
of $X$-filtered $R$-modules where $f$ and $g$ are both bicontrolled maps, the combination of $E$ being lean and $E''$ being insular implies that $E'$ is necessarily scattered.
\end{DefRefName}

\begin{NotRef}{MUDS}
We will use the notation $\frak{C}$ for the class of all coarsely coherent metric spaces.
\end{NotRef}

There is a version of this condition for metric families.  For this definition we use the terminology from Guentner \cite{eG:14}.  We find the notion of the total space of the family especially convenient for our purposes.

\begin{DefRef}{FamGue}
A \textit{metric family} $\{ X_{\alpha} \}$ is simply a collection of metric spaces $X_{\alpha}$.
In all situations in this paper, a metric family will be a collection of subspaces of a given metric space, each equipped with the subspace metric.
The \textit{total space} $X$ of the family $\{ X_{\alpha} \}$ is the disjoint union of the metric spaces $X_{\alpha}$ with the extended metric with values $\infty$ between points from $X_{\alpha}$ and $X_{\beta}$ for $\alpha \ne \beta$.
\end{DefRef}

\begin{DefRefName}{CoarseCoherenceFam}{Coarse Coherence for Families}
A metric family $\{ X_{\alpha} \}$ is \textit{coarsely coherent} if the total space $X$ of the family is coarsely coherent.  

It might be instructive to spell out what this entails. The total space is coarsely coherent if for a collection of exact sequences 
\[
0 \to E'_{\alpha} \xrightarrow{\ f_{\alpha} \ } E_{\alpha} \xrightarrow{\ g_{\alpha} \ } E''_{\alpha} \to 0
\]
of $X_{\alpha}$-filtered $R$-modules where all $E_{\alpha}$ are $D$-lean, all $E''_{\alpha}$ are $d$-insular, all $f_{\alpha}$ and $g_{\alpha}$ are all $b$-bicontrolled maps for some fixed constants $D$, $d$, $b \ge 0$, it follows that all $E'_{\alpha}$ are $\partial$-\textit{scattered} for some uniform constant $\partial \ge 0$.

When we say that a family is in $\frak{C}$, we mean that the family is coarsely coherent.  This is definitely a stronger assumption than each space in the family being in $\frak{C}$.
\end{DefRefName}

The following is the main permanence result of the paper. 

A map between metric spaces $f \colon X \to  Y$ is \textit{uniformly expansive} if there is a function $\phi \colon [0, \infty) \to [0, \infty)$ such that $d_Y (f(x_1), f(x_2)) \le \phi (d_X (x_1, x_2))$ for all pairs of points $x_1$, $x_2$ from $X$.

\begin{ThmRefName}{CCFP}{Fibering Permanence for Coarse Coherence}
Assume $\pi \colon X \to  Y$ is a uniformly expansive map with $Y$ in $\frak{C}$. If for any $r > 0$ the family $\{ f^{-1} (y[r]) \mid y \in Y \}$ is in $\frak{C}$, then $X$ is in $\frak{C}$.
\end{ThmRefName}

The proof of this theorem will require establishing other permanence theorems which we do in the next section.
We will also point out that there is a stronger theorem that can be formulated in terms of families which is a consequence of Theorem \refT{CCFP} and which has many other permanence theorems as corollaries.

\section{Permanence properties and applications}

Let us first introduce additional terminology from coarse geometry.

Two functions $h_1$, $h_2 \colon X \to Y$ between metric space are \textit{close} if there is a constant $C \ge 0$ so that $d_Y (h_1 (x), h_2(x)) \le C$ for all choices of $x$ in $X$.
A function $k \colon X \to Y$ is a \textit{coarse equivalence} if it is uniformly expansive and there exists a uniformly expansive function $l \colon Y \to X$ so that the compositions $k \circ l$ and $l \circ k$ are close to the identity maps.

An example of a coarse equivalence is the notion of quasi-isometry.  This is simply a coarse equivalence $k$ for which the uniformly expansive functions for $k$ and its coarse inverse can be chosen to be linear polynomials.  In geometric group theory, it is very useful that any two choices for a finite generating set of a group produce quasi-isometric word metrics.

The following are basic permanence properties of coarse coherence.

\begin{ThmRefName}{CICC}{Coarse Invariance of Coarse Coherence}
If $X$ and $Y$ are coarsely equivalent then $X$ is coarsely coherent if and only if $Y$ is coarsely coherent.
\end{ThmRefName}

\begin{ThmRefName}{SPCC}{Subspace Permanence for Coarse Coherence}
If $X$ is a subspace of $Y$, and $Y$ is coarsely coherent, then  $X$ is coarsely coherent.
\end{ThmRefName}

From Lemma 6.1 of Guentner \cite{eG:14}, the combination of the two statements is true if and only if a different coarse geometric condition holds.
We recall that a map $k \colon X \to Y$ is a \textit{coarse embedding} if $k$ is a uniformly expansive map which is a coarse equivalence onto its image.  In this situation, the uniformly expansive counterpart $l \colon \im (k) \to X$ is called a \textit{coarse inverse}.
Now the combination of Theorems \refT{CICC} and \refT{SPCC} is true if and only if whenever $Y$ is coarsely coherent and $X$ coarsely embeds in $Y$ then $X$ is coarsely coherent, so it suffices to prove the latter statement.

\begin{proof}
	Suppose $k \colon X \to Y$ is a coarse embedding controlled by the function $\ell$. We assume that $Y$ is coarsely coherent, and we are given an  exact sequence 
\[
0 \to E'_X \xrightarrow{\ f \ } E_X \xrightarrow{\ g \ } E''_X \to 0
\]
of $X$-filtered modules where $f$ and $g$ are both $b$-bicontrolled, $E_X$ is $D$-lean and $E''_X$ is $d$-insular. We aim to show that $E'_X$ is scattered.  Consider the $Y$-filtrations of the modules induced by $k$ as follows: $E_Y (S) = E_X (k^{-1} (S))$.  It follows that $E_Y$ is $\ell (D)$-lean, and $E''_Y$ is $\ell (d)$-insular. Also $f$ and $g$ are bicontrolled as morphisms between $Y$-filtered modules.  We can conclude that $E'_Y$ is $\delta$-scattered for some $\delta \ge 0$, since $Y$ is coarsely coherent.  Now using the same kind of estimate, $E'_X$ is $\ell (\delta)$-scattered.
\end{proof}

There are two types of natural filtrations that can be assigned to submodules of filtered modules.  

\begin{DefRef}{FILTR}
(1) Suppose $F$ is a filtered module and $F'$ is any submodule of $F$. Then $F'$ can be given the canonical filtration $F'(S) = F(S) \cap F'$.  It is clear that if $F$ is an insular filtered module then $F'$ is also insular. 

(2) Suppose $T$ is a subset of $X$.  For any choice of a number $D \ge 0$, one has the submodules $F_{T,D} (S) = \sum_{x \in S \cap T} F(x[D])$ for all $S \subset X$.  They give a filtration  of $F_{T,D} = \sum_{x \in T} F(x[D])$.
It follows easily that $F_{T,D}$ is always a lean $X$-filtered module.  

Notice that, as defined, both filtrations are $X$-filtrations.
\end{DefRef}

We are ready to prove the Fibering Permanence Theorem \refT{CCFP}.

\begin{proof}[Proof of Theorem \refT{CCFP}]
First, observe that in view of Theorem \refT{SPCC} we may assume that $\pi \colon X \to Y$ is surjective.
Given an exact sequence 
\[
0 \to E' \xrightarrow{\ f \ } E \xrightarrow{\ g \ } E'' \to 0
\]
of $X$-filtered $R$-modules where $f$ and $g$ are $b$-bicontrolled maps, $E$ is $D$-lean, and $E''$ is $d$-insular, we want to show $E'$ is scattered.
There is a $Y$-filtration of $E'$ given by $E'_Y (T) = E' (\pi^{-1}(T))$.  It is easy to see that the exact sequence 
\[
0 \to E'_Y \xrightarrow{\ f \ } E_Y \xrightarrow{\ g \ } E''_Y \to 0
\]
of $Y$-filtered modules has the same properties with respect to the new filtrations: $f$ and $g$ are bicontrolled, $E_Y$ is lean, and $E''_Y$ is insular.  This allows to conclude that $E'_Y$ is $\delta_Y$-scattered for some number  $\delta_Y \ge 0$, so every $k \in E'_Y$ is a sum $k = \sum k_y$ where $k_y \in E'_Y (y[\delta_Y]) = E' (\pi^{-1}(y[\delta_Y]))$.

Suppose $\ell$ is a uniform expansion control function for $\pi$.
Let $\mathcal{E}_y$ be the filtered module $E_{\pi^{-1}(y[\delta_Y  + \ell(b)]),D}$, described as option (2) in Definition \refD{FILTR}, which is a lean module containing $f E'_Y (y[\delta_Y])$.
The kernel element $k_y$ is in the kernel of the restriction map $g \colon \mathcal{E}_y \to E''_Y (y[\delta_Y + 2\ell(b) + \ell(D)])$.
The image $\mathcal{E}''_y = g(\mathcal{E}_y)$ is given the canonical filtration induced from the insular filtration of $E''$ which makes $\mathcal{E}''_y$ insular. By the family assumption, we conclude that the kernel of $g \colon \mathcal{E}_y \to \mathcal{E}''_y$ is $\delta$-scattered for a constant $\delta \ge 0$ independent from $y$.  So 
\[
k_y \in \sum_{x \in \pi^{-1}(y[\delta_Y])} E' (x[\delta]) 
\]
because of our assumption that $\pi$ is surjective.
 The conclusion is that $E'$ is $\delta$-scattered.
\end{proof} 

Recall that straight finite decomposition complexity (sFDC) is a property of metric spaces introduced by A. Dranishnikov and M. Zarichnyi \cite{aDmZ:12}. The class of groups with sFDC is remarkably broad.  It includes groups with finite asymptotic dimension, all elementary amenable groups, and all countable subgroups of almost connected Lie groups.

\begin{ThmRef}{KBZAAA}
	A metric space $X$ with sFDC is coarsely coherent.
\end{ThmRef}

\begin{proof}
	The proof of the main result from \cite{bG:13} applies literally and verifies coarse coherence of the metric space.
\end{proof}

It is known that a metric space with finite asymptotic dimension has sFDC.

\begin{CorRef}{bkdkfj}
	A metric space $X$ with finite asymptotic dimension is coarsely coherent.
\end{CorRef}

The following simple corollary to Theorem \refT{CCFP} shows that coarse coherence of finitely generated groups is preserved by group extensions. 

\begin{CorRef}{OUR}
Let $\pi \colon G \to H$ be a surjective homomorphism from a finitely generated group $G$.  We assume that the groups are given word metrics with respect to finite generating sets and that the kernel $K$ is given the subspace metric.
If $K$ is coarsely coherent and $H$ is coarsely coherent then $G$ is coarsely coherent.
\end{CorRef}

\begin{proof}
If $S$ is a finite generating set for $G$, $\pi (S)$ can be used as a finite generating set of $H$, and the resulting word metric space is known to be coarsely coherent by quasi-isometry invariance of coarse coherence.  The isometric action of $G$ on $H$ is transported from the left action of $G$ on itself, so the action is by isometries.  In this situation, the quasi-stabilizers $W_r (e)$ from the proof of \cite[Theorem 7]{gBaD:06} are $\pi^{-1}(e[r]) = K[r]$.  Since $K$ is coarsely equivalent to $K[r]$, we have $\pi^{-1}(e[r])$ is coarsely coherent for any $r$ and so all $\pi^{-1}(x[r])$ are coarsely coherent by quasi-isometry invariance.  This means that for a fixed $r>0$ the family $\{ \pi^{-1}(x[R]) \}_{x \in X}$ is uniformly coarsely coherent.  We conclude that $G$ is coarsely coherent by applying the Fibering Permanence Theorem \refT{CCFP}.
\end{proof}

Using more of the recent technology developed by D. Kasprowski, A. Nicas, and D. Rosenthal in \cite{dKaNdR:18} we can deduce many other permanence properties for coarse coherence.  As a consequence, coarse coherence is preserved under many other group theoretic constructions.  

First we observe that Theorem \refT{SPCC} together with the Subspace Permanence from Theorem \refT{SPCC} confirm the following General Fibering Permanence property of coarse coherence.  For the sake of brevity we don't review the extension of uniformly expansive maps of metric families.  The details can be found in \cite{dKaNdR:18}.

\begin{ThmRef}{CCFPnn}
Assume $\pi \colon \mathcal{X} \to \mathcal{Y}$ is a uniformly expansive map with $\mathcal{Y}$ in $\frak{C}$. If for any uniformly bounded subfamily $\mathcal{B}$ of $\mathcal{Y}$ we have that the family $f^{-1} (\mathcal{B})$ is necessarily in $\frak{C}$, then $X$ is in $\frak{C}$.
\end{ThmRef}

\begin{RemRef}{JK}
	The remarkable sequence of Theorems 5.4, 5.6, 5.8, 5.10 and 5.12 from \cite{dKaNdR:18} establishes a number of permanence properties as formal consequences of Theorem \refT{CCFPnn} for a collection of metric families $\frak{C}$ which also contains all metric families  with finite asymptotic dimension.  We have seen that this is so for coarse coherence from Corollary \refC{bkdkfj}.  The new properties that follow are the Finite Amalgamation Permanence, Finite Union Permanence, Union Permanence and Limit Permanence.  We refer to \cite{dKaNdR:18} for the precise definitions and details.
\end{RemRef}

\section{Coarse regular coherence}

We want to define the new property \textit{coarse regular coherence} of finitely generated groups and leverage our permanence theorems for coarse coherence to say what we know about the class of coarsely regular coherent groups.
We also want to recall the notion of \textit{weak regular coherence} from \cite{gCbG:04,gCbG:16,bG:13} and explain the relationship to  this paper.

In the rest of this section, $R$ will be a commutative Noetherian ring.

Let $\Gamma$ be a finitely generated group with a word metric.  There is an isometric action of $\Gamma$ on itself by left multiplication.

\begin{DefRef}{HUT}
A $\Gamma$-filtered $R$-module is $\Gamma$-\textit{equivariant} or simply a $\Gamma$-\textit{module} if $F(\gamma S) = \gamma
F(S)$ for all choices of $\gamma \in \Gamma$ and $S \subset \Gamma$.
\end{DefRef}

It is clear that a $\Gamma$-module has the structure of an $R[\Gamma]$-module.    
It turns out that any $R[\Gamma]$-module can be given a $\Gamma$-filtration specific to a finite generating set in $\Gamma$ as follows.
Given a left $R[\Gamma]$-module $F$ with a finite generating set
$\Sigma$, it is also an $R$-module with the generating set
${\Sigma}' = \{ \gamma \sigma \in F \mid \gamma \in \Gamma, \sigma \in \Sigma \}$.
Now one can associate to every subset $S$ of $\Gamma$
the left $R$-submodule $F(S)$ generated by $\gamma \sigma \in {\Sigma}'$ such that
$\gamma \in S$ and $\sigma \in \Sigma$.
This gives a functor $F \colon \mathcal{P}(\Gamma) \to \mathrm{Mod}_R (F)$, from
the power set of $\Gamma$ to the $R$-submodules of $F$ such that $F(\Gamma) = F$, $F(\emptyset) = 0$, and for a bounded subset
$T \subset \Gamma$, $F(T)$ is a finitely generated $R$-module.
This shows that $F$ with the given filtration, which will be denoted by $F_{\Sigma}$, is a $\Gamma$-filtered $R$-module.
It is easy to see that $F_{\Sigma}$ is $\Gamma$-equivariant, so $F_{\Sigma}$ is a $\Gamma$-module.
Clearly, $F_{\Sigma}$ is $0$-lean by design.

As an example, a finitely generated free $R[\Gamma]$-module with a finite set of free generators $\Sigma$ can be given the structure of a $\Gamma$-module as above.  In this case ${\Sigma}' = \Gamma \times \Sigma$.
It is easy to see that in this case the $\Gamma$-filtration is $0$-lean and $0$-insular.
It follows from Theorem \refT{lninpres} that images of idempotents of lean, insular, locally finitely generated modules also have these three properties.  This allows us to generate examples of such modules from idempotents of free geometric modules.
Examples of harder to construct non-projective lean, insular $R[\Gamma]$-modules can be found in the last section of \cite{gCbG:16}.

\medskip

Coarse coherence of the group $\Gamma$ implies that the kernel of any $R[\Gamma]$-equivariant surjection between finitely generated, lean, insular $\Gamma$-modules is finitely generated.  We refer to \cite{gCbG:04} for an explanation.  So we have the following consequence.

\begin{ThmRef}{resWC}
Suppose $R$ is a commutative Noetherian ring and $\Gamma$ is a finitely generated group viewed as a metric space with respect to the word metric associated to some fixed choice of a generating set. If $\Gamma$ is coarsely coherent then any finitely generated $\Gamma$-module $F$ has a resolution of finite type by finitely generated free $\Gamma$-modules.
\end{ThmRef}

\begin{proof}
Given $F$, select a finite generating set $\Sigma$ and consider the $\Gamma$-filtered module $F_{\Sigma}$ as above.  Now the $R$-submodule $\langle \Sigma \rangle$ generated by $\Sigma$ has a free finitely generated $R$-module $F_e$ with a surjection onto $\langle \Sigma \rangle$.  Similarly, all submodules $ \langle \gamma \Sigma \rangle$ have free finitely generated $R$-modules $F_{\gamma}$ isomorphic to $F_e$ and surjections onto each $\langle \gamma \Sigma \rangle$ which serve as components of a $\Gamma$-equivariant surjection 
$\bigoplus_{\gamma \in S} F_{\gamma} \to F$.  Coarse coherence of $\Gamma$ gives that the kernel of this surjection is finitely generated and is, in fact, lean when given the canonical filtration as a submodule of the free geometric $\Gamma$-module $\bigoplus F_{\gamma}$.  This allows to iterate the construction inductively to produce a resolution of $F$ of finite type.
\end{proof}

Now here is finally the definition.

\begin{DefRefName}{CRCoherence}{Coarse Regular Coherence}
A finitely generated group $\Gamma$ is \textit{coarsely regular coherent relative to a ring $R$} if every $R[\Gamma]$-module $F$ which is lean, insular, and locally finitely generated when viewed as a $\Gamma$-filtered $R$-module has a finite projective resolution over $R[\Gamma]$.  It is called simply \textit{coarsely regular coherent} if it is coarsely regular coherent relative to any regular coherent ring of finite global dimension.
\end{DefRefName}

\begin{NotRef}{MUDSGHT}
We will use the notation $\frak{Y}$ for the class of all coarsely regular coherent groups.
\end{NotRef}

Now let us address the relationship of coarse coherence to weak coherence.

\begin{DefRef}{LODG}
A $\Gamma$-module $F$ has an \textit{admissible presentation} if there is an exact sequence 
\[
F_2 \to F_1 \to F \to 0
\]
where $F_1$ and $F_2$ are finitely generated free $\Gamma$-modules, and the homomorphism $f \colon F_2 \to F_1$ is bicontrolled.
\end{DefRef}

\begin{DefRefName}{jkiu1}{Weak Regular Coherence}
A finitely generated group $\Gamma$ is \textit{weakly coherent relative to $R$} if every $\Gamma$-module with an admissible presentation has a projective resolution of finite type over $R[\Gamma]$. 
It is \textit{weakly regular coherent relative to $R$} if every $\Gamma$-module with an admissible presentation has a finite projective resolution over $R[\Gamma]$.  
We define $\Gamma$ to be simply \textit{weakly coherent}, resp. \textit{weakly regular coherent}, if $\Gamma$ is weakly coherent, resp. weakly regular coherent, relative to any regular Noetherian ring of finite global dimension.
\end{DefRefName}

Weak regular coherence was introduced and studied in \cite{gCbG:04a,gCbG:04,gCbG:16,bG:13}.

\begin{ThmRef}{MVXEDF}
	Coarse coherence of a group implies weak coherence. Coarse regular coherence implies weak regular coherence. 
\end{ThmRef}

The proof follows from the next two lemmas.

\begin{LemRef}{LeanBC2}
A module with an admissible presentation is lean, insular, and locally finitely generated with respect to the $\Gamma$-filtration.
If $\Gamma$ is a coarsely coherent group then conversely,
every lean, insular, finitely generated
$R[\Gamma]$-module
has an admissible presentation.
\end{LemRef}

\begin{proof}[Proof of Lemma \refL{LeanBC2}]
If $F$ is the quotient of a boundedly bicontrolled homomorphism $F_2 \to F_1$ with the image $I$ filtered by $I(S) = I \cap F_1 (S)$, then the filtration $F(S) = F_1 (S)/I (S)$ makes the quotient map $F_1 \to F$ bicontrolled.  It is known that boundedly bicontrolled homomorphisms are balanced (Proposition 2.8 of \cite{gCbG:11}), so the kernel $K$ of this quotient map is lean as the image of a lean module $F_2$.  Therefore $F$ is lean and insular by parts (1) and (3) of Theorem \refT{lninpres} applied to the short exact sequence $K \to F_1 \to F$.

Now suppose $F$ is a $D$-lean, insular, finitely generated $R[\Gamma]$-module.
Since the $R$-submodule $F(e[D])$ is finitely generated, there is a finitely generated free $R$-module $F_e$ and an $R$-linear epimorphism $\phi_e \colon F_e \to F(e[D])$.  One similarly has an epimorphisms $\phi_{\gamma} \colon F_{\gamma} \to F(\gamma[D])$ using isomorphic copies $F_{\gamma}$ of $F_e$.
A new geometric $\Gamma$-filtered module $F_1$ is defined by assigning $F_1 (S) = \bigoplus_{\gamma \in S} F_{\gamma}$, so $F_1$ is a lean insular $\Gamma$-module.
There is an $R[\Gamma]$-homomorphism $\phi_1 \colon F_1 \to F$ induced by sending $F_1 (\gamma)$ onto $F(\gamma[D])$ via $\phi_{\gamma}$.
This is a $D$-bicontrolled $R$-homomorphism. 
Since $\Gamma$ is coarsely coherent, the kernel $K$ of $\phi_1$ is a scattered $\Gamma$-module.
It is insular by part (2) of Theorem \refT{lninpres}.
One similarly constructs a filtered module $F_2$ and a bicontrolled $R[\Gamma]$-homomorphism $\phi_2 \colon F_2 \to K$ using the $\Gamma$-equivariant scattering of $K$.
Since $K$ is finitely generated over $R[\Gamma]$, $F_2$ can be chosen to be finitely generated as well.
The composition of $\phi_2$ and the inclusion of $K$ in $F_1$ gives a bicontrolled $\psi_2 \colon F_2 \to F_1$ with the cokernel $F$, as required.
\end{proof}

\begin{LemRef}{LeanBC3}
If $\Gamma$ be a coarsely coherent group then
every lean, insular, finitely generated
$R[\Gamma]$-module
is of type $\mathrm{FP}_{\infty}$.
\end{LemRef}

\begin{proof}[Proof of Lemma \refL{LeanBC3}]
Applying the construction from the proof of Lemma \refL{LeanBC2} to $\phi_2 \colon F_2 \to K$ and proceeding inductively,
one obtains a projective $R[\Gamma]$-resolution by finitely generated $R[\Gamma]$-modules
\[
\ldots \rightarrow F_{n} \rightarrow F_{n-1}
\rightarrow \ldots \rightarrow F_2 \rightarrow F_1 \rightarrow F
\rightarrow 0
\]
In this resolution, each homomorphism $\psi_n \colon F_n \to F_{n-1}$
factors through an epimorphism $\phi_n$ onto a scattered, insular $\Gamma$-submodule $K_{n-1}$ of $F_{n-1}$.
\end{proof}

To state the next theorem which provides examples of coarsely regular coherent groups, recall that P.~Kropholler \cite{pK:93,pK:99} defined the class of groups $\mathrm{LH} \mathcal{F}$ which includes, in particular, all groups with finite $K(\Gamma,1)$.

\begin{ThmRef}{KLORSCMK}
A coarsely coherent group that belongs to Kropholler's hierarchy $\mathrm{LH} \mathcal{F}$ is coarsely regular coherent.
\end{ThmRef}

\begin{proof}
	The main result of \cite{bG:13} was stated for weak regular coherence. From inspecting the argument, it is evident that the key property of the syzygies analyzed in the proof is that they are scattered rather than lean.  So the argument applies verbatim under this new assumption.
\end{proof}

\begin{CorRef}{sFDCyes}
A finitely generated group with sFDC that belongs to Kropholler's hierarchy $\mathrm{LH} \mathcal{F}$ is coarsely regular coherent.
\end{CorRef}

Finally, we list some closure properties of the class $\frak{Y}$.

\begin{ThmRef}{HAWLKI}
The class of coarsely regular coherent groups is closed under
\begin{itemize} 
\item passage to subgroups,
\item passage to supergroups of finite index,
\item extensions such as finite semi-direct products of groups, including finite direct products,  
\item direct unions, 
\item amalgamated free products and HNN extensions.
\end{itemize}
\end{ThmRef}

\begin{proof}
	   The corresponding closure properties for $\mathrm{LH} \mathcal{F}$ are in sections 2.2 and 2.4 of \cite{pK:93}.  For the class $\frak{C}$ with the coarse coherence property, these are consequences of Theorem \refT{KLORSCMK} in conjunction with Theorem \refT{CCFPnn} and Corollary \refC{bkdkfj} as explained in Remark \refR{JK}.
\end{proof}

Note an elementary consequence of this fact: the family of regular coherent groups that Waldhausen constructed in \cite{fW:78}, which are essentially multiple amalgamated free products and HNN extensions of finitely many copies of $\mathbb{Z}$, are in $\frak{Y}$.

\end{document}